\documentclass{amsart}

\usepackage{amsthm}
\usepackage{graphicx}
\usepackage{amssymb}

\newtheorem{teo}{Theorem}[section]
\newtheorem{lemma}[teo]{Lemma}
\newtheorem{prop}[teo]{Proposition}
\newtheorem{cor}[teo]{Corollary}

\theoremstyle{definition}
\newtheorem{defi}[teo]{Definition}
\newtheorem{rem}[teo]{Remark}

\theoremstyle{remark}
\newtheorem{prof}[teo]{Proof of}

\def\mr{\mathbb{R}}
\def\mz{\mathbb{Z}}

\begin{document}
\title{Branched shadows and complex structures on $4$-manifolds}

\author[Costantino]{Francesco Costantino}
\address{Institut de Recherche Math\'ematique Avanc\'ee\\
  7, Rue Ren\'e Descartes\\
  67084 Strasbourg Cedex France}
\email{f.costantino@sns.it}


\begin{abstract}
We define and study branched shadows of $4$-manifolds as a combination of branched spines of $3$-manifolds and of Turaev's shadows. We use these objects to combinatorially represent $4$-manifolds equipped with $Spin^c$-structures and homotopy classes of almost complex structures. We then use branched shadows to study complex $4$-manifolds and prove that each almost complex structure on a $4$-dimensional handlebody is homotopic to a complex one.
\end{abstract}

\maketitle

\tableofcontents
\section{Introduction}
Shadows were defined by V. Turaev at the
beginning of the nineties in \cite{Tu3} as a method for
representing knots alternative to the standard one based on knot
diagrams and Reidemeister moves. The theory was developed in
the preprint ``Topology of shadows", later included as
a revised version in \cite{Tu}; a short account of the
theory was also published in \cite{Tu2}. Since then, only few applications of shadows
were studied. Among these we recall the use of
shadows made by U.
Burri in \cite{Bu} and A. Shumakovitch in \cite{Sh} 
to study Jones-Vassiliev invariants of knots and the study
of ``Interdependent modifications of links and invariants of
finite degree" developed by M.N. Goussarov in \cite{Gou}.
More recently, in \cite{CT} and \cite{CT2} D.P. Thurston and the author use shadows to study $3$-manifolds, their geometry and to prove a quadratic upper bound on the minimal complexity of a $4$-manifold whose boundary is a given $3$-manifold.

On the other side, {\it branched polyhedra} (and in particular branched spines of $3$-manifolds) were studied in depth by Benedetti and Petronio since 1997 in \cite{BP}. They showed how branched spines encode $3$-manifolds equipped with additional topological structures as $Spin^c$-structures, vector fields and a more refined structure they called {\it concave traversing fields}. In \cite{BPco}, they provided a calculus for these objects, i.e. a finite set of local modifications connecting any two branched spines of the same manifold respecting the extra-structure. Later, in \cite{BPto}, using these objects they extended the definition of the Reidemeister-Turaev refined torsion to $3$-manifolds with arbitrary boundary. Branched spines were also used by Benedetti and Baseilhac to define quantum hyperbolic invariants in \cite{BB} and \cite{BB2}.

The present paper is devoted to combine the notions of shadow of a $4$-manifold and branched polyhedron: the resulting theory is a generalization in dimension $4$ of a set of results which hold for branched spines of $3$-manifolds. The main interest of the new theory is that a series of these results can be further developed and studied in new directions which in dimension $3$ are not visible as, for instance, the study of the relations between branched shadows and complex structures.

After an introductory section in which we recall the notion of shadow and how to thicken it to a $4$-manifold, we show how a branching encodes a $Spin^c$-structure on that manifold. Then we prove in Theorem \ref{teo:surjectivityofrefinedreconstruction} a result that we summarize as follows:
\begin{teo}
Let $M$ be a $4$-dimensional handlebody (i.e. a $4$-manifold admitting a handle decomposition without $3$ and $4$-handles). Each $Spin^c$-structure and each homotopy class of almost complex structure on $M$ can be encoded by a branched shadow of $M$.  
\end{teo} 

Hence branched shadows are a key tool for a combinatorial approach to a series of $4$-dimensional problems related with $Spin^c$-structures and almost complex structures. 
We further develop the theory in the last section by studying complex structures and branched shadows. It turns out that a branched shadow embedded in a complex manifold behaves like an embedded oriented real surface. 
We translate to the world of shadows a series of classical definitions and results due to Chern, Spanier, Bishop and Lai on indices of embedded real surfaces. We prove a shadow version of the well known result of Harlamov and Eliashberg which allows, under suitable conditions, to annihilate pairs of complex points in embedded real surfaces. 
In the end, as an application, we prove the following:
 \begin{teo}\label{mainteo}
 Let $M$ be a $4$-dimensional handlebody. Each homotopy class of almost complex structures on $M$ can be represented by an integrable complex structure.
 \end{teo}
Our proof is self-contained and based on the machinery set up in the present paper; the same result was proved by P. Landweber in \cite{La} through different techniques.
In a subsequent paper (\cite{Co4}) we will study sufficient combinatorial conditions assuring that a branched shadow determines a Stein domain.

{\bf Acknowledgements.} The author wishes to warmly thank Stephane Baseilhac, Riccardo
Benedetti, Paolo Lisca, Dylan Thurston and Vladimir Turaev for their criticism
and encouraging comments.

\section{Shadows of $4$-manifolds}

In this section we recall the notion of shadow of a $4$-manifold;
for a complete account of this topic we refer to \cite{Tu} and, for an introductory one
\cite{Co5}. From now on, all the manifolds and homeomorphisms will be smooth unless explicitly stated. 

\begin{figure}
  \centerline{\includegraphics[width=8.4cm]{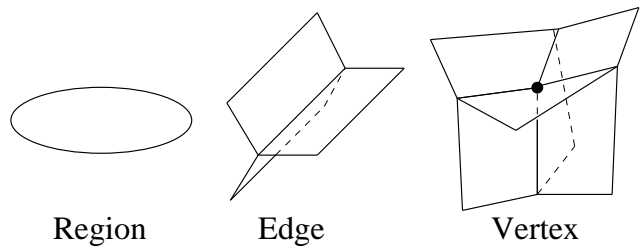}}
  \caption{The three local models of a simple polyhedron. }
  \label{fig:singularityinspine}
\end{figure}

\begin{defi} A simple polyhedron $P$ is a $2$-dimensional CW complex whose
local models are those depicted in Figure
\ref{fig:singularityinspine}; the set of points whose neighborhoods have models
of the two rightmost types is a $4$-valent graph, called
{\it singular set} of the polyhedron and denoted by $Sing(P)$. The
connected components of $P-Sing(P)$ are the {\it regions} of $P$.
A simple polyhedron whose singular set is connected
and whose regions are all discs is called {\it
standard}. 
\end{defi}
From now on, for the sake of simplicity all the
polyhedra will be standard without
explicit mentioning.

\begin{defi}[Shadow of a $4$-manifold]\label{def:flat}
  Let $M$ be a smooth, compact and oriented $4$-manifold.  $P \subset N$ is a
  \emph{shadow} for $M$ if:
  \begin{enumerate}
  \item $P$ is a closed polyhedron  embedded in $M$ so that $M-P$ is diffeomorphic to $\partial M\times (0,1]$;
  \item $P$ is \emph{flat} in $M$, that is for each point
  $p\in P$ there exists a local chart $(U,\phi)$ of $M$ around $p$ such
  that $\phi(P\cap U)$ is contained in $\mr^3 \subset \mr^4$ and in this chart the
pair $(\mathbb{R}^3\cap \phi (U),\mathbb{R}^3\cap \phi(U\cap P))$ is diffeomorphic
to one of the models depicted in Figure \ref{fig:singularityinspine}.
\end{enumerate}
\end{defi}

From now on, all the embedded polyhedra will be flat unless explicitly stated.
\begin{rem}
Note that the original definition of shadows was given in the PL setting by Turaev in \cite{Tu}
. But in four dimensions the smooth and the PL setting are equivalent, that is for each PL-structure on a compact manifold there exists a unique compatible (in a suitable sense) smooth structure. Definition \ref{def:flat} is the natural translation of the notion of shadow to the smooth setting. 
\end{rem}

A necessary and sufficient
condition (see \cite{Co}) for a $4$-manifold $M$ to admit a shadow is that $M$ is a $4$-handlebody, that is $M$ admits a handle decomposition without $3$ and $4$-handles.
In particular, $\partial M$ is a non empty connected $3$-manifold.
From now on, all the manifolds will be $4$-handlebodies unless explicitly stated. 

\begin{figure}
  \centerline{\includegraphics[width=8.4cm]{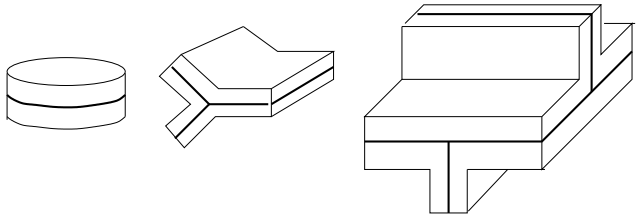}}
  \caption{The three type of blocks used to thicken a spine of a $3$-manifold.  }
  \label{fig:spineblocks}
\end{figure}

Can we reconstruct the neighborhoods of a polyhedron $P$ in a manifold $M$ from its
combinatorics? Let us first understand the easier $3$-dimensional
case, where the polyhedron is embedded in an oriented
$3$-manifold $N$ and its combinatorics
allows one to reconstruct its
regular neighborhoods in the following way. Any decomposition of 
$P$ in the local
patterns of Figure \ref{fig:singularityinspine}, induces a
decomposition of any of its regular neighborhoods in blocks as those of
Figure \ref{fig:spineblocks}. These can be reglued to each
other according to the combinatorics of $P$. That way, a
polyhedron embedded in a $3$-manifold, determines
its regular neighborhoods in the $3$-manifold. It is known \cite{Cas} that
any $3$-manifold with non-empty boundary can be reconstructed that
way, as a neighborhood of some embedded standard polyhedra, called {\it spine}.

Let us now pass to the $4$-dimensional case. Suppose that $P$ is a surface embedded in a oriented
$4$-manifold $M$. In general we cannot reconstruct the tubular
neighborhoods of $P$ by using only its topology, since their structure depends on the self-intersection number of $P$ in $M$.  To state it
differently, the tubular neighborhood of a surface in a
$4$-manifold is homeomorphic to the total space of a disc bundle
over the surface (its normal bundle), and the Euler number of this
bundle is a necessary datum to reconstruct its topology.

Hence, we see that to encode the topology of a
neighborhood of $P$ in $M$ we need to``"decorate" $P$ with
some additional information; when $P$ is an oriented surface,
the Euler number of its normal bundle is a sufficient datum.

We describe now the basic decorations we need for a general standard polyhedron $P$, for a more detailed account see \cite{Tu}. 
Let us denote $\mathbb{Z}[\frac{1}{2}]$ the group of integer multiples of $\frac{1}{2}$. There are two canonical {\it colorings} on the regions of $P$, i.e. assignments of elements of $\mathbb{Z}_2$ or $\mathbb{Z}[\frac{1}{2}]$, the second depending on a flat embedding of $P$ in an oriented  $4$-manifold. They are:

{\bf The $\mathbb{Z}_2$-gleam of $P$}, constructed as follows.  Let $D$ be the (open) $2$-cell associated to a given region of $P$ and $\overline{D}$ be
the natural compactification $\overline{D}=D\cup S^1$ of the (open) surface represented by
$D$. The embedding of $D$ in $P$ extends to a map
$i:\overline{D}\to P$ which is injective in $int(\overline{D})$,
locally injective on $\partial \overline{D}$ and which sends
$\partial \overline{D}$ into $Sing(P)$. Using the map $i$ we can
``pull back" a small open neighborhood of $D$ in $P$ and
construct a simple polyhedron $N(D)$ collapsing on $\overline{D}$
and such that the map $i$ extends as a local homeomorphism
$i':N(D)\to P$ whose image is contained in a small
neighborhood of the closure of $D$ in $P$. When $i$ is an embedding of $\overline{D}$ in $P$, then $N(D)$
turns out to be homeomorphic to a neighborhood of $D$ in $P$
and $i'$ is its embedding in $P$. In general, $N(D)$ has the following structure: each boundary
component of $\overline{D}$ is glued to the core of an annulus or
of a M\"obius strip and some small discs are glued along half of
their boundary on segments which are properly embedded in these
annuli or strips and cut transversally once their cores. We
define the $\mathbb{Z}_2$-gleam of $D$ in $P$ as the
reduction modulo $2$ of the number of M\"obius strips used to
construct $N(D)$. This coloring only depends on
the combinatorial structure of $P$.

\begin{figure}
  \centerline{\includegraphics[width=6.4cm]{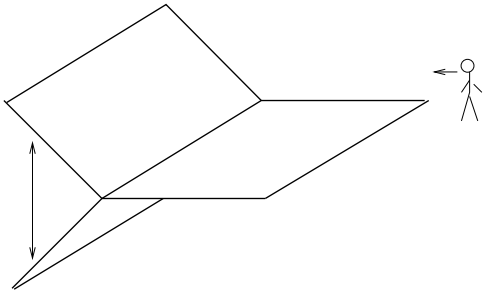}}
  \caption{The picture sketches the position of the polyhedron in a
    $3$-dimensional slice of the ambient $4$-manifold.  The direction
    indicated by the vertical double arrow is the one along which the
    two regions touching the horizontal one get separated.  }
  \label{fig:divergingdirection}
\end{figure}

{\bf The gleam of $P$}\label{embpoly}, constructed as follows. Let us now suppose that $P$ is flat in an oriented $4$-manifold $M$, with $D$, $\overline{D}$ and $i:\overline{D}\rightarrow P$ as above. 
Pulling back through $i$ a small neighborhood of $i(N(D))$ in $M$, we obtain a $4$-dimensional neighborhood of $N(D)$. This is an oriented ball $B$ on which we fix an auxiliary riemannian metric. Since $N(D)$ is locally flat in $B$, $N(D) - D$ well defines a line normal to $\overline{D}$ in $B$ along $\partial \overline{D}$ and hence a section of the projectivized normal bundle of $\overline{D}$ (see Figure \ref{fig:divergingdirection}). Let then $gl(D)$ be equal to $\frac{1}{2}$ times the obstruction to extend this section to the whole $\overline{D}$; such an obstruction is an element of $H^2(\overline{D},\partial \overline{D};\pi_1(S^1))$, which is canonically identified with $\mathbb{Z}$ since $B$ is oriented. Note that the gleam of a region is integer if and only if its $\mathbb{Z}_2$-gleam is zero. 

Using the fact
that the $\mathbb{Z}_2$-gleam is always defined, Turaev 
generalized~\cite{Tu} the
notion of gleam to non-embedded polyhedra as follows:
\begin{defi}
A {\it gleam} on a simple polyhedron $P$ is a coloring on the regions of
$P$ with values in $\mathbb{Z}[\frac{1}{2}]$ such that the color of
a region is integer if and only if its $\mathbb{Z}_2$-gleam is
zero.
\end{defi}
\begin{teo}[Reconstruction Theorem \cite{Tu}]\label{teo:reconstruction}
Let $P$ be a polyhedron with gleams~$gl$; there exists a
canonical reconstruction map associating to $(P,gl)$ a pair $(M_P,P)$
where $M_P$ is a smooth, compact and oriented $4$-manifold, and $P\subset M$ is a shadow of $M$ (see Definition \ref{def:flat}).
If $P$ is a standard polyhedron flat in a smooth and oriented  $4$-manifold and $gl$ is the gleam of $P$ induced by its embedding, then $M_P$ is diffeomorphic to a compact neighborhood of $P$ in $M$.
\end{teo}

The proof is
based on a block by block reconstruction procedure similar to the
one used to describe $3$-manifolds by means of their spines.
Namely, for each of the three local patterns of Figure \ref{fig:singularityinspine}, we consider the $4$-dimensional thickening given by the
product of an interval with the corresponding $3$-dimensional
block shown in Figure \ref{fig:spineblocks}. All
these thickenings are glued to each other according to the
combinatorics of $P$ and its gleam.

By Theorem \ref{teo:reconstruction}, to study $4$-manifolds one can either
use abstract polyhedra equipped with
gleams or embedded polyhedra. The latter approach is more
abstract, while the former is purely combinatorial; 
we will use both approaches in
the following sections. The translation in the combinatorial setting of the definition of shadow of a $4$-manifold is the following:

\begin{defi}[Combinatorial shadow]
A polyhedron equipped with gleams $(P,gl)$ is said to be a {\it
shadow} of the $4$-manifold $M$ if $M$ is diffeomorphic to the manifold associated to $(P,gl)$ by means of
the reconstruction map of Theorem
\ref{teo:reconstruction}.
\end{defi}

Given a shadow $(P,gl)$ of a $4$-manifold $M$, it is possible to
modify it by a series of local modifications called ``moves".
\begin{figure}[htbp]
  \centerline{\includegraphics[width=11.4cm]{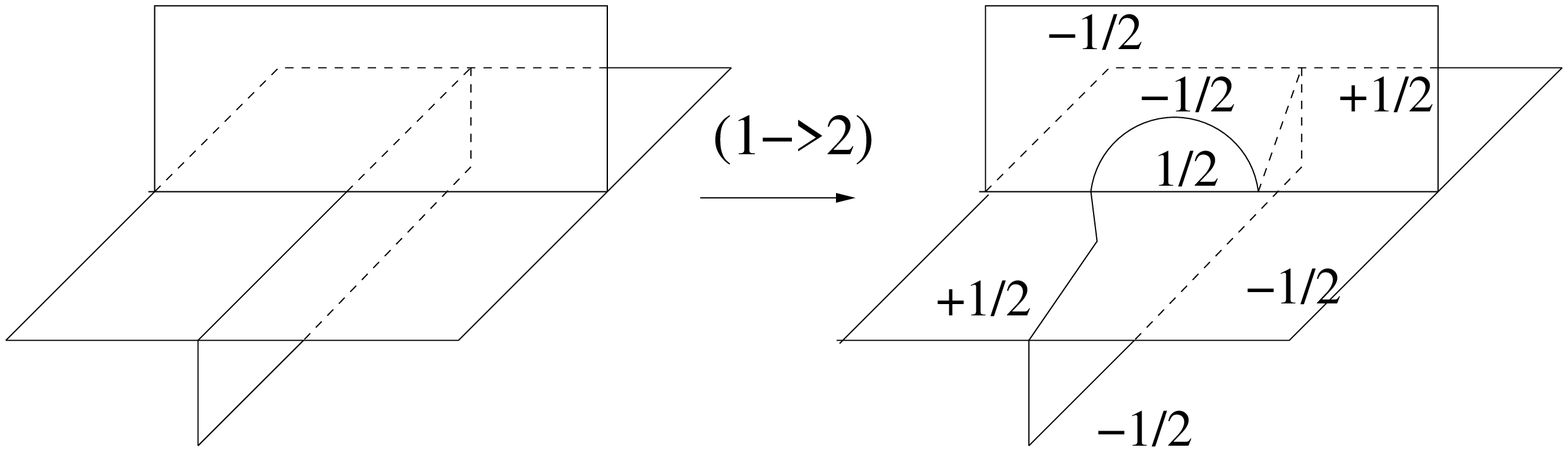}}
  \centerline{\includegraphics[width=11.4cm]{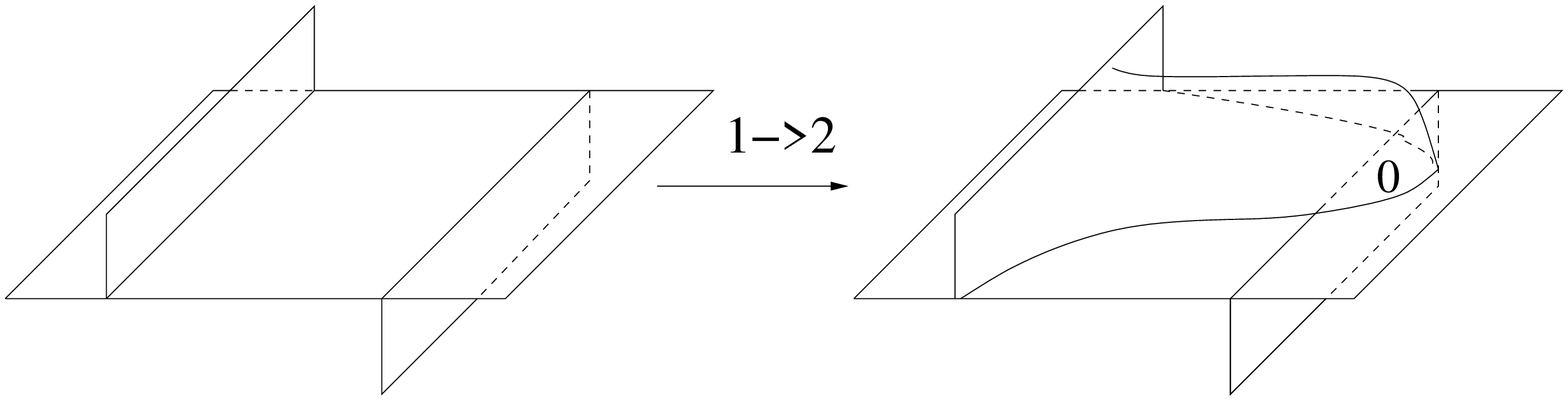}}
  \centerline{\includegraphics[width=11.4cm]{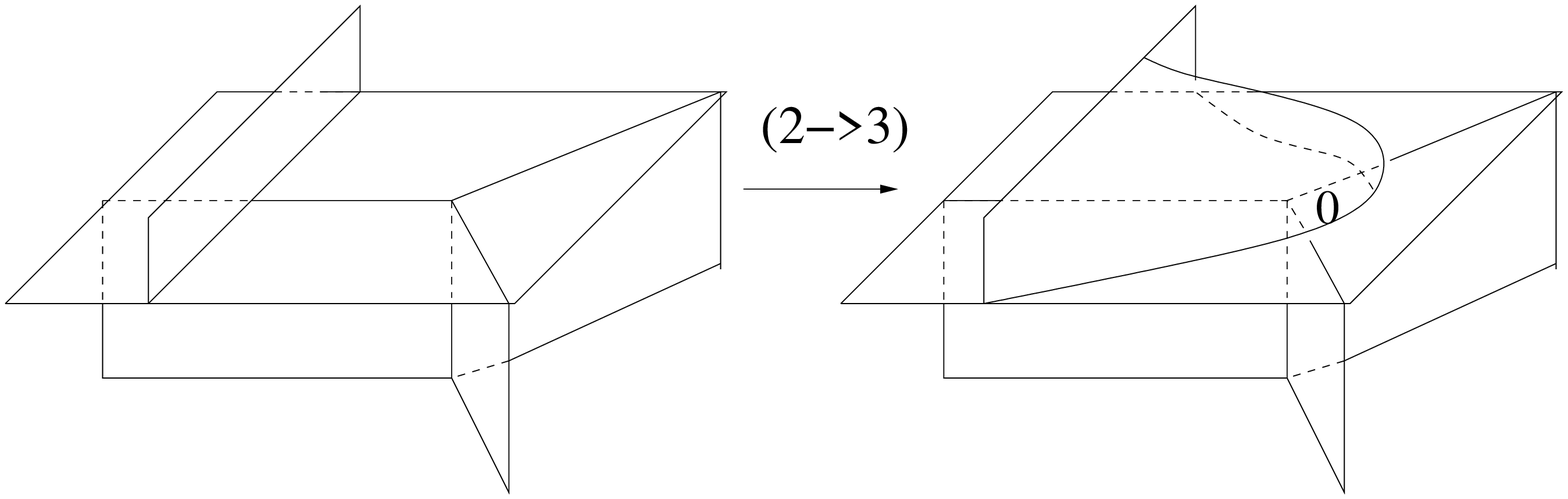}}
  \caption{The three shadow equivalences.}
  \label{fig:2-3}
\end{figure}
The three most important are shown in Figure
\ref{fig:2-3}.
To visualize these moves, imagine that a region of $P$
slides over some other regions, producing a polyhedron which differs from the initial one only in a
contractible subpolyhedron (drawn in the left part of
the figure). These moves are called respectively $1\to 2$, $0\to
2$ (or lune or finger-move) and $2\to 3$-moves (or Matveev-Piergallini move), because of their effect on the
number of vertices of the polyhedra. 

Let us analyze first
the $0\to 2$-move; this
move acts in a $4$-ball contained in $M$ and containing the part
of $P$ shown in the left part of the picture. After the move, we
replace this part of $P$ with that drawn in the
right and obtain a shadow $P'$ of $M$. The whole move can be performed in a $3$-dimensional slice of the $4$-ball in $M$.  
The same
comments apply to the case of the $2\to
3$-move.

The case of the $1\to 2$-move is slightly different; this
move starts applies to a neighborhood of a vertex and slides a region (the
vertical lower one in the left part of the figure) over the
vertex thus creating a new one. It is a good exercise to
visualize this sliding: it cannot be performed in $\mr^3$. 

The above comments apply also to the inverses of the
moves. We now clarify the meaning of the
numbers written on the regions of the polyhedra in the figure. Each move represents a modification of
the embedded polyhedron, and produces a new polyhedron whose gleam
is induced by the embedding in the ambient manifold as explained in the definition of the gleam. Each region of this polyhedron corresponds to one of the initial polyhedron except the disc created by the positive moves and contained in the right part of the pictures.
The numbers in the figure represent the change in the gleams of
the regions; the gleam of the new disc is explicited in the pictures. The gleam of the regions which do not appear in the figure do not change since their embedding in $M$ does not.

All these comments can be summarized in the following proposition:
\begin{prop}
Let $(P,gl)$ be a shadow of a $4$-manifold $M$ and let $(P',gl')$
be obtained from $(P,gl)$ through the application of any sequence
of $1\to 2$, $0\to 2$ and $2\to 3$ moves and their inverses. 
Then $(P',gl')$ is a shadow of $M$. 
\end{prop}
Note that $(P',gl')$ can be reconstructed from $(P,gl)$ by following step by step the changes induced by the moves of the sequence.

\section{Branched shadows}
 Given a simple polyhedron $P$ we define the notion of
{\it branching} on it as follows:
\begin{defi}[Branching condition]\label{branchingcondition}
A branching $b$ on $P$ is a choice of an orientation for each
region of $P$ such that for each edge of $P$, the
orientations induced on the edge by the regions containing it do
not coincide.
\end{defi}

\begin{rem}
This definition corresponds to the definition of ``orientable
branching'' of \cite{BP}.
\end{rem}

We say that a polyhedron is {\it branchable} if it admits a
branching and we call {\it branched polyhedron} a pair
$(P,b)$ where $b$ is a branching on $P$.
\begin{defi}
Let $(P,gl)$ be a shadow of a $4$-manifold $M$. We
call {\it branched shadow} of $M$ the triple $(P,gl,b)$ where
$(P,gl)$ is a shadow and $b$ is a branching on $P$. 
When this will not cause any confusion, we will not
specify the branching $b$ and we will simply write $(P,gl)$.
\end{defi}
\begin{prop}\label{prop:existsbranchedshadow}
Any $4$-manifold admitting a shadow admits also a branched shadow.
\end{prop}
\begin{prf}{1}{
We sketch the idea of the proof which is an adaptation of Theorem 3.4.9 of \cite{BP}. We note that a branched shadow $P'$ is obtained from a shadow $P$ via an algorithmic procedure.
Orient arbitrarily all the regions of $P$; if for each edge of
$Sing(P)$ the three orientations induced on it by the 
regions containing it do not coincide, then we already found a
branching on $P$. Let us suppose then that an edge $e$ is induced
three times the same orientation,
and that the endpoints of $e$ lie on two different
vertices of $P$ (we can always find a shadow for which all the
edges have this property by applying a suitable sequence of 
the moves of Figure \ref{fig:2-3} on $(P,gl)$). The basic idea of
the proof is to apply a $2\to 3$-move along the edge $e$ to ``blow
up" it and create a new region whose orientation we can choose
arbitrarily: indeed a $2\to 3$-move makes an edge disappear and
creates a new region whose boundary is formed by three edges (see Figure \ref{fig:2-3}).
Choosing appropriately the orientation of the new region and keeping unchanged the orientations of the other ones, we
can assure that no edge is induced three times the same orientation
from the regions containing it. In some particular cases, one of
the edges touched by the new region does not satisfy the
branching condition, hence we apply again a $2\to
3$-move on it. In \cite{BP} an accurate analysis of the possible cases is performed
showing that this process eventually ends with a branched
polyhedron. }
\end{prf}

\begin{figure} [h!]
   \centerline{\includegraphics[width=11.4cm]{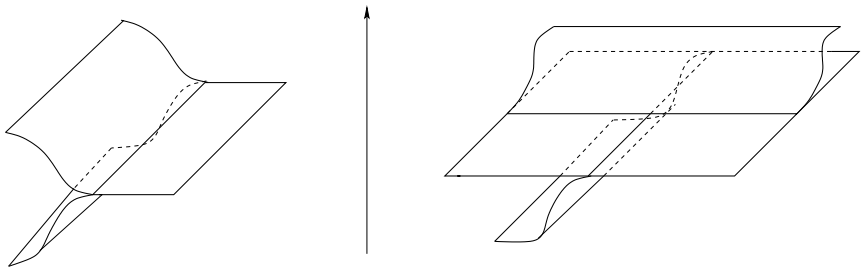}}
  \caption{How a branching allows a smoothing of the polyhedron: the
  regions are oriented so that their projection on the ``horizontal" (orthogonal to the drawn vertical direction) plane is orientation preserving.}\label{branching}
\end{figure}

A branching on a shadow allows us to smoothen its singularities
and equip it with a smooth structure as shown in Figure
\ref{branching}. This smoothing can be performed also inside the
ambient manifold obtained by thickening the shadow; the shadow locally appears as in 
Figure \ref{branching}, where the two regions orienting the edge
in the same direction approach each other so that, for any auxiliary riemannian metric on the ambient manifold, all the derivatives of their
distance go to zero while approaching the edge.

If $(P,gl)$ is a branched shadow of a $4$-manifold, and we apply to it one of the moves of Figure \ref{fig:2-3}, we get another
shadow $P'$ of the same manifold containing one region more than $P$. 
Each region
of $P$ naturally corresponds to a region of $P'$ and the
region of $P'$ which does not correspond to one of $P$ is the small 
disc created by the move (see Figure
\ref{fig:2-3}). Hence the branching on $P$ induces a choice of an orientation
on each region
of $P'$ except on that disc and these orientations satisfy the branching condition
on all the edges of $P'$ not touched by that disc. Analogously, if $P'$ is
obtained from $P$ through the inverse of a basic move, then each region of
$P'$ corresponds to a region of $P$ and the branching on $P$
induces an orientation on each region of $P'$.

\begin{defi}\label{branchedmove}
A basic move $P\rightarrow P'$ applied on a branched polyhedron $P$ is called {\it
  branchable} if
it is possible to choose an orientation on the disc created by the move which,
together with the orientations on the regions of $P'$ induced by the branching of $P$,
defines a branching on $P'$. Analogously, the inverse of a basic move
  applied to $P$ is branchable if the orientations induced by the
  branching of $P$ on the regions of $P'$ define a branching.
\end{defi}

\begin{figure} [h!]
   \centerline{\includegraphics[width=8.4cm]{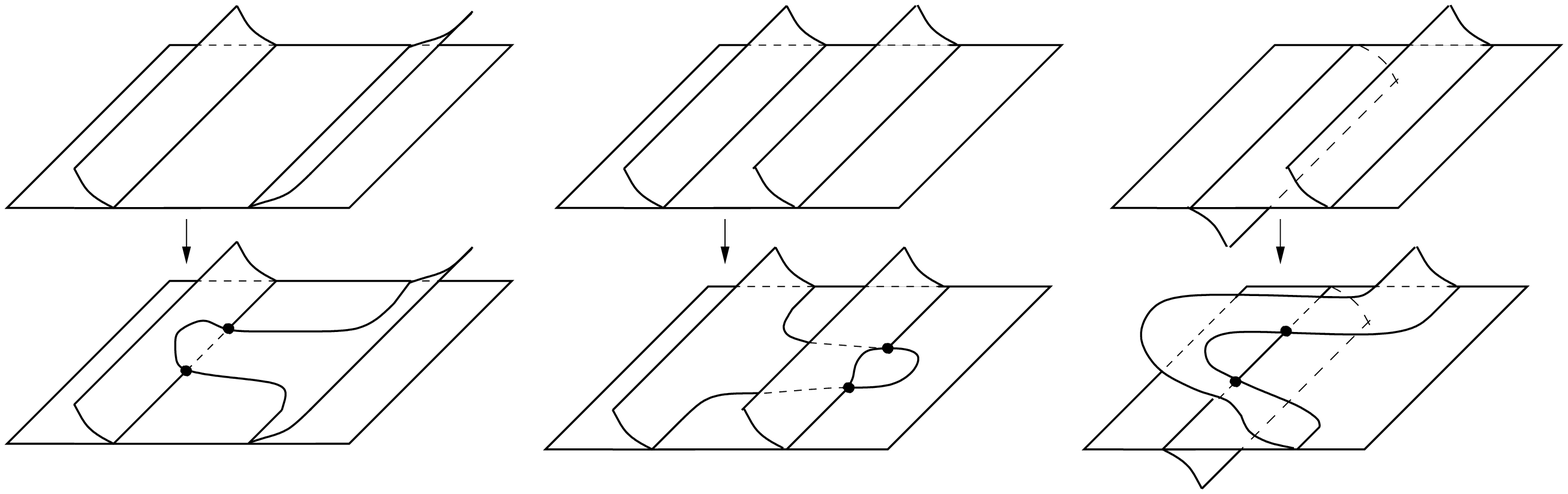}
\includegraphics[width=5.4cm]{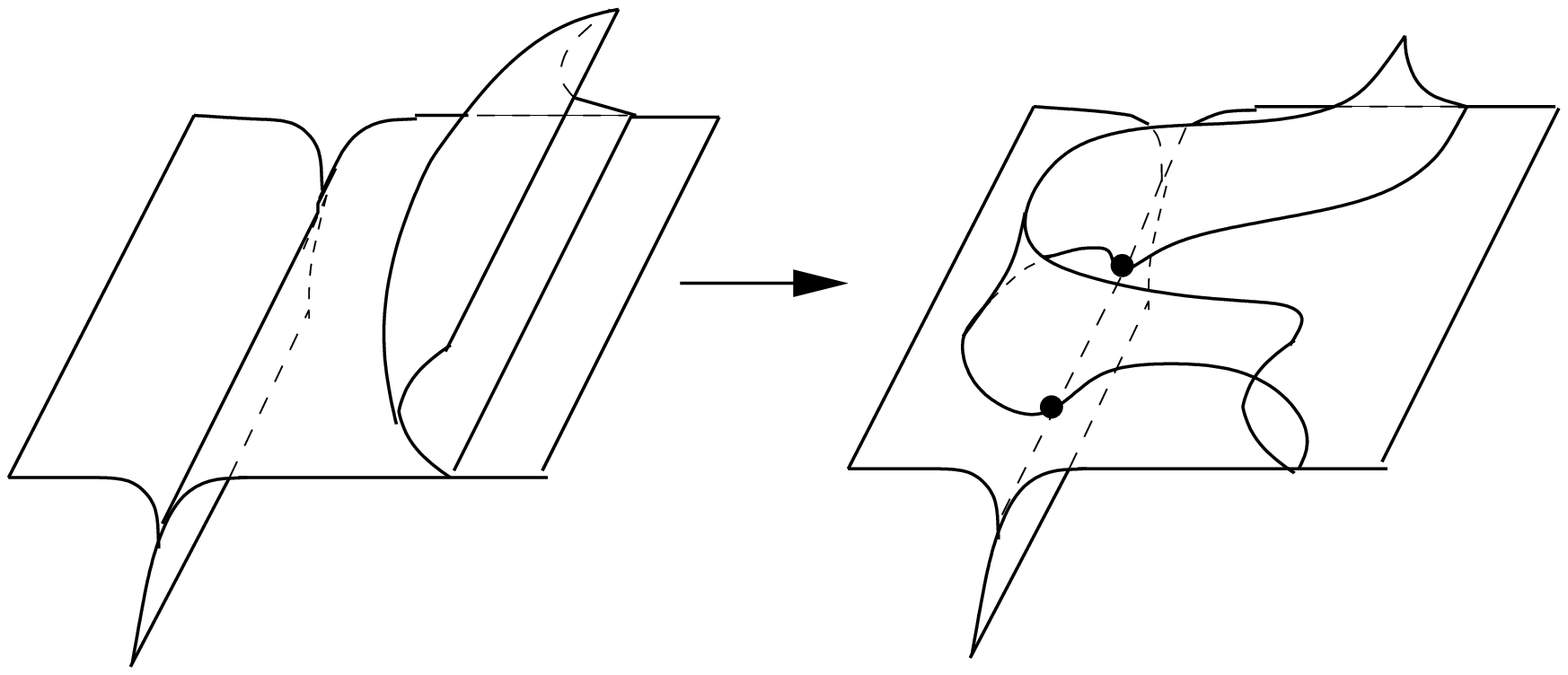}}
  \caption{In the left part of this figure we show the three branched
  versions of the lune-move called
  ``sliding''-moves. In the right part we show the version
  called ``bumping''-move.}\label{branched lune}
\end{figure}

\begin{figure} [h!]
   \centerline{\includegraphics[width=11.4cm]{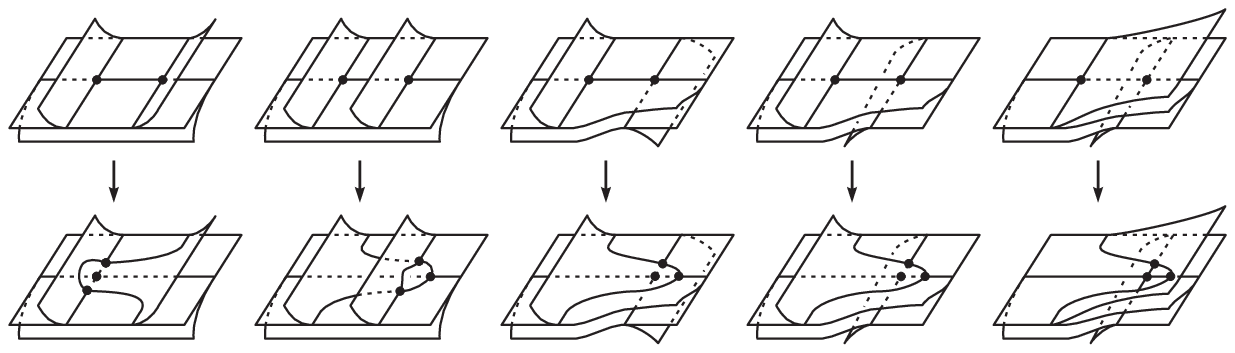}}
   \centerline{\includegraphics[width=8.4cm]{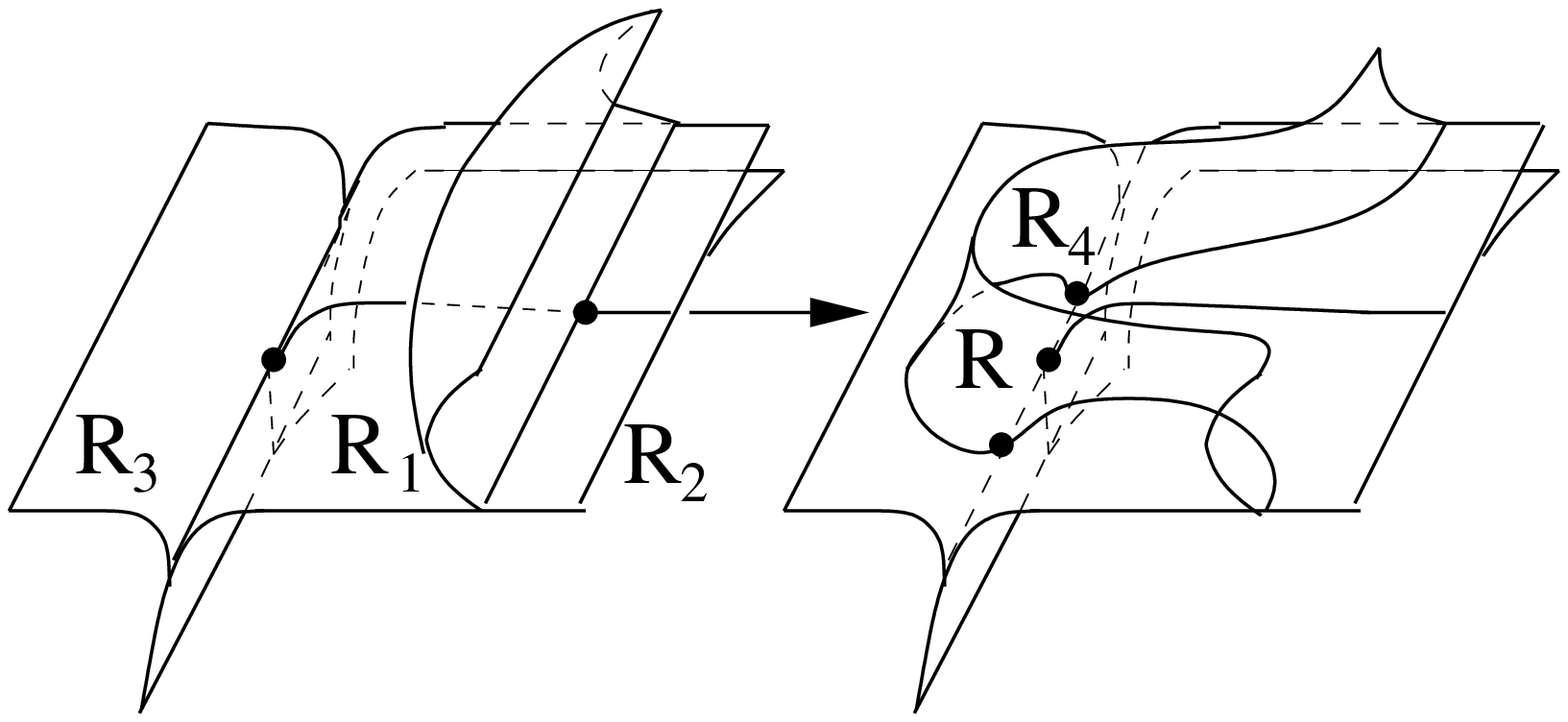}}

  \caption{In the upper part of the figure we show the $5$
 branched versions of the $2\to 3$-move called
``sliding''-moves. The bottom part of the figure represents the ``bumping''-move.
 }\label{branched MP}
\end{figure}

A branching is a kind of loss of symmetry on a
polyhedron and this is reflected by the fact that each move has
many different branched versions.
To enumerate all the possible branched versions of the moves, one has
to fix any possible
orientation on the regions of the left part of Figure
\ref{fig:2-3} and complete these
orientations in the right part of the figure by fixing one orientation 
on the region created by the move; by Definition \ref{branchedmove},
one obtains a branched version of a basic move when the branching
condition is satisfied both in the left and in the right part of
the figure.
Fortunately, many of the possible combinations are equivalent
up to symmetries of the pictures;
we show in Figure \ref{branched lune} all the branched versions of
the $0\to 2$-move and in Figure \ref{branched MP} those of the $2\to 3$-move.
In these
figures we split these branched versions in two types namely the
{\it sliding}-moves and the {\it bumping}-moves; this
differentiation will be used later.

We will also use often the following terminology:
\begin{defi}
Let $e$ be an edge of a branched polyhedron $P$ and let $R_i$,$R_j$ and $R_k$ be the regions of $P$ containing it in their boundary. Then $R_i$ is said to be the {\it preferred region of $e$} if it induces the opposite orientation on $e$ with respect to those induced by $R_j$ and $R_k$.
\end{defi}

\begin{figure}
   \centerline{\includegraphics[width=6.4cm]{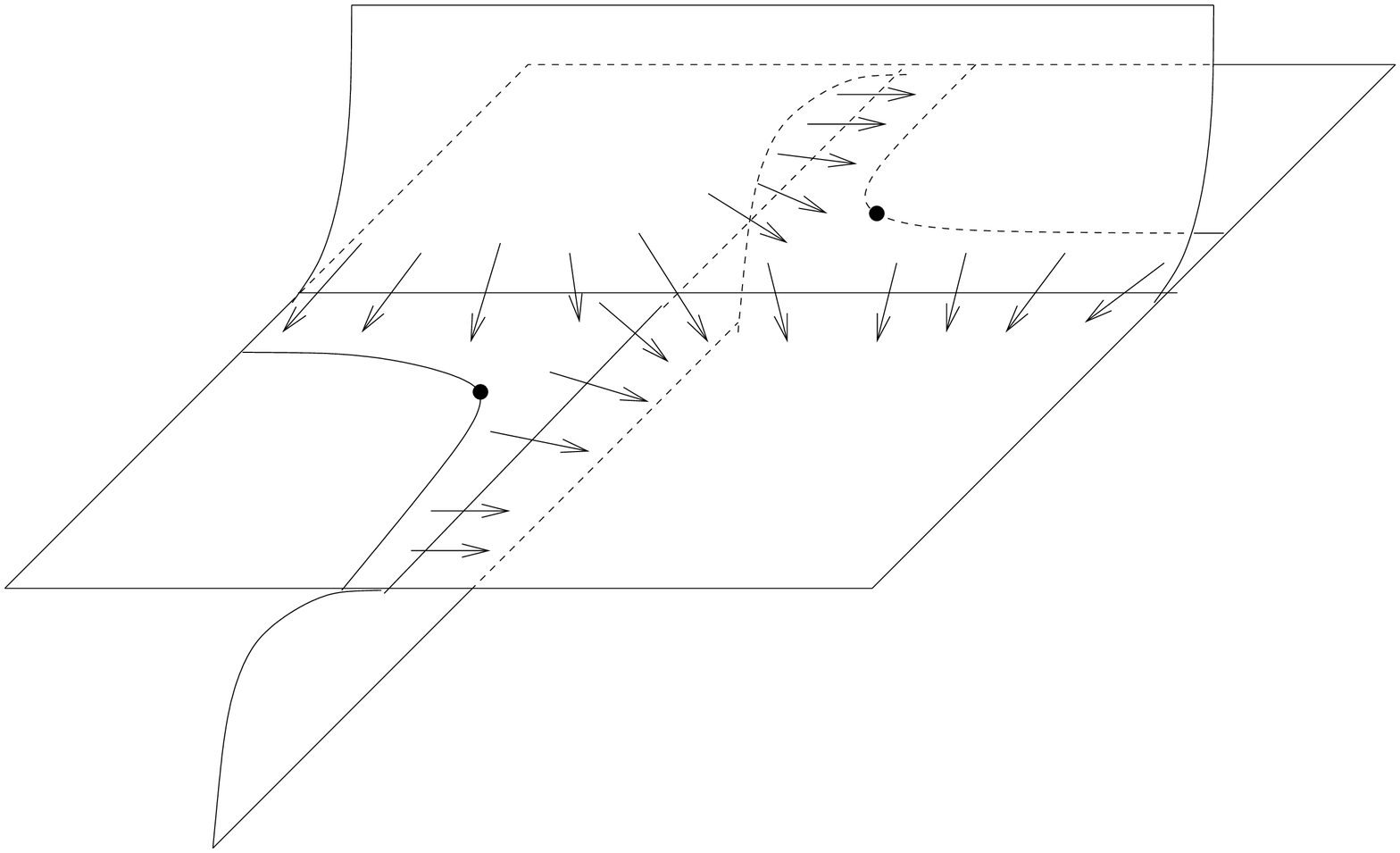}}

  \caption{ In this figure we show how the maw behaves near a
  vertex of a branched polyhedron.
  }\label{fig:maw}
\end{figure}

It has been proved in \cite{BP}[Chapter 3] that any lune
and $2\to 3$-move is branchable, but some of their inverses are not. 
Regarding the $1\to 2$-move and its inverse, the
following holds:
\begin{lemma}\label{lem:p1branchable}
Each $1\to 2$-move or its reverse is branchable.
\end{lemma}
\begin{prf}{1}{
\begin{figure} [h!]
   \centerline{\includegraphics[width=12.4cm]{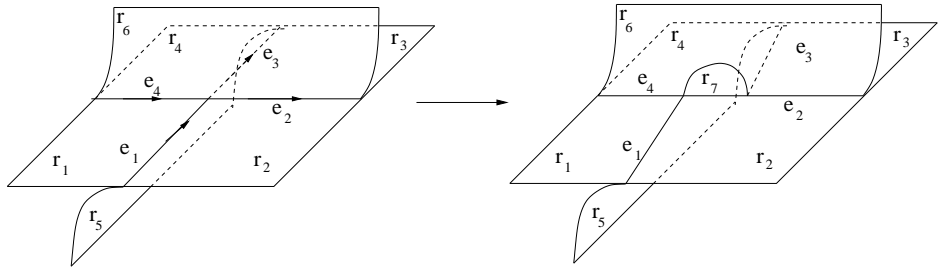}}

  \caption{In this picture we fix the notation we use in the proof of Lemma \ref{lem:p1branchable}.}\label{p1branchable}
\end{figure}
Let $P$ be a branched polyhedron on which a $1\to 2$-move acts. To fix the notation, we use the right-hand rule in Figure \ref{p1branchable} and give a
sign to the each region of the polyhedron involved in the $1\to 2
$-move: $+$ if
the region is positively oriented with respect to the upward
direction and $-$ otherwise. That way, we identify a branching
near the vertex with a six-uple of signs.
There are 24 possibilities but since the opposite of
a branching is a branching we reduce to study
the 12 branchings where the sign of the region $r_1$
is $+$. A branching induces an
orientation on the edges $e_i,\ 1=1,2,3,4$ touching the vertex;
in the case represented in the left part of the figure all the signs are $+$.

We need
to show that, for each branching (and then for each
six-uple of signs) near the vertex in the left part of Figure
\ref{p1branchable}, there is a choice of an orientation (and so a
sign) for the region $r_7$ in the right part of the figure so that
no edge $e_i,\ i=1,2,3,4,5,6$ is induced three times
the same orientation from the regions containing it.
This is proved in the following table where, in the left column, we list the
possible branchings before the move, and in the right column we
write the signs corresponding to the compatible orientations of $r_7$.
Note that there are cases where both the
orientations of $r_7$ give a branching to $P$.
\\

\begin{tabular}{ll}

\begin{tabular}{|l|l|}
\hline $(r_1,r_2,r_3,r_4,r_5,r_6)$ & $r_7$\\ \hline
$(+,+,+,+,+,+)$ & $+$\\ $(+,+,+,+,+,-)$ & $\pm$ \\ $(+,+,+,+,-,+)$
& $\pm$ \\ $(+,+,+,+,-,-)$ & $-$ \\ $(+,+,+,-,+,+)$ & $+$ \\
$(+,+,-,-,+,+)$ & $+$ \\ \hline
\end{tabular}

&

\begin{tabular}{|l|l|}
\hline $(r_1,r_2,r_3,r_4,r_5,r_6)$ & $r_7$\\ \hline
$(+,+,-,+,-,+)$ & $\pm$\\ $(+,-,+,+,-,-)$ & $-$ \\ $(+,+,-,-,-,+)$
& $+$ \\ $(+,-,-,+,-,-)$ & $-$ \\ $(+,-,-,+,-,+)$ & $-$ \\
$(+,-,-,-,-,+)$ & $\pm$ \\ \hline
\end{tabular}
\end{tabular}
\\

The case of the $(1\to 2) ^{-1}$-move is simpler: if the polyhedron is branched before the
inverse move, then it is after, since no new edge is created
during the move and no regions merge.
 }
\end{prf}
The following proposition is a consequence of the fact that a $4$-handlebody retracts on its shadows.
\begin{prop} \label{prop:branched presentation}
Let $P$ be a branched shadow of a $4$-manifold
$M$, and let $R_i,\ i=1,\ldots,n$ and $e_j,j=1,\ldots,m$ be
respectively the regions and the edges of $P$ oriented according
to the branching of $P$. Then $H_2(M;\mz)$ is the kernel of
the boundary application $\partial :\mz[R_1,\ldots,R_n] \to
\mz[e_1,\ldots,e_m]$. Moreover
$H^{2}(M;\mz)$ is the abelian group generated by
the cochains $\hat{R}_i,\ i=1,\ldots,n$ dual to the regions of $P$
subject to the relations generated by the coboundaries of the edges having the form  $\delta(\hat{e}_j)=-\hat{R}_i+\hat{R}_j+\hat{R}_k$ where $R_i$ is the preferred region of $e_j$. \end{prop}
Given a shadow $(P,gl)$ of a $4$-manifold $M$, there are three cochains representing classes in $H^2(M;\mathbb{Z})$ naturally associated to $(P,gl)$. 
The first one is the {\it Euler cochain} of $P$, denoted $Eul(P)$ and constructed as follows.
Let $m$ be the vector field tangent to $P$ (using the smoothed structure given by the branching)
which near the center of the edges points inside the preferred
regions; we extend $m$ in a neighborhood of the
vertices as shown in Figure \ref{fig:maw}.

The field constructed above is the {\it maw}
of $P$. For each region $R_i$ of $P$, the maw
gives a vector field defined near $\partial R_i$, so
it is possible to extend this field to a tangent field on
the whole $R_i$ having isolated singularities of indices
$\pm 1$; let $n_i$ be the algebraic sum of these indices over the
region $R_i$. The Euler cochain is defined as $Eul(P)=\sum_i
n_i\hat{R}_i$; its meaning will be analyzed in the next section.

\section{Branched shadows and almost complex structures}\label{sec4}

As before, let $M$ be an oriented $4$-handlebody and let $g$ be a fixed auxiliary riemannian metric on $M$. In this section we show that a
branched shadow determines a pair $(M,[J])$ where
$[J]$ is a homotopy class of almost complex structures on $M$ suitably compatible with $g$,
and that for each such class, there exists a branched shadow of
$M$ encoding it. 

\begin{defi}
An almost complex structure $J$ on an oriented $4$-manifold $M$ is a smooth morphism $J:TM\to TM$ such that for each point $p\in M$ it holds $J^2=-Id$. We say that $J$ is {\it positive} if at each point $p\in M$ there exists a positive (with respect to the orientation of $M$) basis of $T_pM$ of the form $(x,J(x),y,J(y))$. We say that $J$ is orthogonal with respect to $g$ if for each $p\in M$ the map $J:T_pM\to T_pM$ is $g$-orthogonal.
\end{defi}

Let $P$ be a branched shadow of $M$. In each of the local blocks used to reconstruct
$M$ from $P$ as in Theorem \ref{teo:reconstruction} (these blocks are the products of an interval with the
$3$-dimensional blocks of Figure \ref{fig:spineblocks}), $P$ is smoothly embedded in a non symmetrical way (see Figure \ref{branching}). As in Figure \ref{branching}, we 
choose an
{\it horizontal} $2$-plane and
orient it according to the branching of $P$, so that each of
the basic building blocks of $M$ is equipped with a
distribution of oriented $2$-planes, denoted as $T(P)$.

Let $V(P)$ be the field of
oriented {\it vertical} $2$-planes of $TM$ which are pointwise positively 
$g$-orthogonal to the planes of $T(P)$.

We define an almost complex structure $J_P$ by requiring that its restriction to the two fields $T(P)$ and $V(P)$ acts pointwise as a $\frac{\pi}{2}$ positive rotation and we extend this action linearly to the whole $TM$. By construction $J_P$ is positive and $g$-orthogonal. Note that the choice of $g$ was arbitrary. The following lemma assures us that our constructions are well defined up to homotopy. 

\begin{lemma}\label{almostcomplex=orthogonal}
Let $M$ be an oriented $4$-handlebody and $g$ a riemannian metric on $M$. 
If $P$ is a branched shadow of $M$ and $g'$ is another riemannian metric on $M$, the almost complex structures on $J_P$ and $J'_P$ constructed as above using respectively $g$ and $g'$ are homotopic. 
Moreover, each homotopy class of almost complex structures on $M$ contains a $g$-orthogonal representative.
\end{lemma}
\begin{prf}{1}{
Let $J_0$ be an arbitrary almost complex structure on $M$. 
We first prove the second statement and split its proof in two steps:
\begin{enumerate}
\item there exists a riemannian metric $g_0$ on $M$ such
  that $J_0$ is $g_0$ orthogonal;
\item $J_0$ is homotopic to a $g$-orthogonal almost complex
  structure $J_1$. 
\end{enumerate}

{\bf Step 1.} The space of the scalar
metrics on $\mathbb{R}^4$ with respect to which a fixed complex
structure $J$ on $\mathbb{R}^4$ is orthogonal is a non-empty convex
set: indeed in a fixed basis the orthogonality condition can be written
as $J^tgJ$ where $J$ is the complex structure and $g$ is the
symmetric matrix representing the scalar product. In particular, the
fiber bundle over $M$ of pointwise scalar products with 
respect to which $J_0$ is pointwise
orthogonal has contractible
fiber, and hence admits a section $g_0$.

{\bf Step 2.} The homotopy of metrics on $M$ given by
$g_t=(1-t)g_0+tg$ connects $g_0$ to $g$: we show how to
lift it to a homotopy $J_t$ of $g_t$-orthogonal almost
complex structures connecting $J_0$ to a
$g$-orthogonal almost complex structure. 

Consider the bundle of $J_0$ complex
$2$-planes on $M$; its fiber is $\mathbb{CP}^1=S^2$ and so, since $M$
retracts onto a $2$-dimensional polyhedron, we
can find a section and choose a
field $F$ of $J_0$-complex tangent $2$-planes on $M$. 
For each
point $p\in M$ we define $J_t$ to be the $\frac{\pi}{2}$
rotation with respect to the metric $g_t$ both on $F_p$ and on its
$g_t$-orthogonal $2$-plane. This gives us a
$g_t$-orthogonal almost complex structure $J_t$ connecting $J_0$ to a $g$-orthogonal almost complex structure concluding the proof of the second statement.
To prove the first statement, use Step 2 with $F=T(P)$ and $g_0=g'$.
}\end{prf}

\begin{cor}\label{refinedreconstruction}
Let $M$ be a smooth $4$-manifold admitting a shadow. The restriction to the branched shadows of $M$ of the reconstruction map of Theorem \ref{teo:reconstruction} 
 can be refined to a map whose image
 is contained in the set of pairs $(M,[J])$ where $[J]$ is a homotopy class
 of positive almost complex structures on $M$.
\end{cor}

From now on, we will only use positive and $g$-orthogonal almost complex structures. 
Our construction splits the tangent bundle of $M$ as the sum of two
linear complex bundles $V(P)$ and $T(P)$, hence the first Chern
class of $TM$, viewed as a complex bundle using the almost complex
structure $J_P$, is equal to
$c_1(T(P))+c_1(V(P))$. The following proposition
(whose proof is identical to that of Proposition 7.1.1 of \cite{BP}) is a recipe to calculate $c_1(T(P))$:
\begin{prop}\label{prop:eulercochain2}
The class in $H^2(P;\mz)$ represented by the Euler cochain $Eul(P)$
coincides with the first Chern class of the horizontal plane field
$T(P)$ of $P$ in $M$.
\end{prop}

We define the {\it gleam cochain} $gl(P)$ as
$gl(P)=\sum_i gl(R_i)\hat{R}_i$, where $\hat{R}_i$ is the cochain
dual to $R_i$, the coefficient $gl(R_i)$ is the gleam of the
region $R_i$ and the sum ranges over all the regions $R_i$ of $P$.
Note that since the gleams can be half-integers, it is not a
priori obvious that the gleam cochain represents an integer class
in $H^2(M;\mathbb{Z})\cong H^2(P;\mathbb{Z})$. 
\begin{lemma}\label{glclass}
The gleam cochain $gl(P)$ represents in $H^2(M;\mathbb{Z})$ the
first Chern class of the field of oriented $2$-planes $V(P)$.
\end{lemma}
\begin{prf}{1}{
The normal bundle of each region of $P$ is given by
$V(P)$. Call $PV(P)$ the projectivisation of $V(P)$ and let $v$ be a
generic section of $V(P)$ on $P$. The number of
zeros of the projection of $v$ in $PV(P)$ on a $2$-cycle is 
twice the number of zeros of $v$ on the same cycle; indeed, the index of each zero of the projection of $v$ in
$PV(P)$ is the double of the index of the corresponding
zero of $v$. By the construction of the gleam of a region, $2$-times
the gleam cochain represents the first Chern
class of $PV(P)$. In particular, the number of zeros of the
projection of $v$ in $PV(P)$ on a $2$-cycle is $2$ times
the evaluation of the gleam cochain on the cycle. Hence the number
of zeros of $v$ on the same cycle equals the evaluation of the
gleam cochain on the cycle. }
\end{prf}

\begin{cor}\label{Chern}
The first Chern class of the almost complex structure on $M$ associated
with $P$ is represented by the cochain $c_1(P)=Eul(P)+gl(P)$.
\end{cor}

\subsection{Comparing different almost complex structures}

Corollary \ref{refinedreconstruction} does not tell us which homotopy classes $[J]$ of almost complex structures can be
reconstructed by means of suitable branched shadows of $M$: we
prove in the subsequent subsections that the reconstruction
map is surjective on the pairs $(M,[J])$. To do this, we need
to ``compare" different homotopy classes of almost complex
structures on $M$. This can be performed through a shortcut based
on the theory of $Spin^c$-structures.
\begin{teo}[\cite{KM}]\label{teo:bij}
Let $M$ be a smooth, compact and oriented $4$-manifold admitting a
shadow and equipped with an auxiliary orthogonal riemannian metric
$g$. 
\begin{enumerate}
\item There exists a canonical bijection $b$ between the set
$J$ of homotopy classes of orthogonal, positive, almost complex
structures on $M$ and the set $S$ of $Spin^c$-structures on $M$. 
\item Fixed an arbitrary $Spin^c$-structure $s$ on $M$, any other such structure
$s'$ is isomorphic to $s\otimes l$ where $l$ is a complex line bundle
over $M$; this equips the set of $Spin^c$-structures (and hence of
homotopy classes of almost complex structures) on $M$ with a structure of affine space over $H^2(M;\mathbb{Z})$ where $s'-s=c_1(l)$. 
\item To each $Spin^c$-structure $s$ one can associate an element of $H^2(M;\mathbb{Z})$ called its  {\it Chern class} and denoted by $c_1(s)$ such that if $s$ and $s'$ are as above then $c_1(s')-c_1(s)=2c_1(l)$ and such that $c_1(b(J))=c_1(J)$ for each almost complex structure $J$. 
\end{enumerate}
\end{teo}
\begin{cor}
If $H^2(M;\mathbb{Z})$ has no $2$-torsion, then two almost complex structures are homotopic if and only if their Chern classes coincide.
\end{cor}
\begin{prf}{2}{
All the results on $Spin^c$-structures are standard, see \cite{KM}. In that paper a construction is exhibited which associates to each pair $(s,\psi)$ where $s$ is a $Spin^c$-structure on $M$ and $\psi$ a positive spinor on it, an almost complex structure defined out of the zeros of $\psi$. If $\psi$ is generic the almost complex structure is defined out of a finite set. 

If $M$ collapses onto a $2$-dimensional polyhedron, the spinor can be choosen to be non-zero and each two such spinors can be connected through a path of non-zero spinors. This gives the map $b^{-1}$ from $Spin^c$-structures into the set of almost complex structures; it is a standard fact that it is surjective since each almost complex structure $J$ allows one to explicitly construct the $Spin^c$-structure $b(J)$; also the last statement comes directly from the construction of $b(J)$ from $J$.
}\end{prf}

In particular, a branched shadow $P$ can be
used to reconstruct a $Spin^c$-structure on $M$ instead of the
homotopy class $[J_P]$. A direct consequence of Theorem \ref{teo:bij} 
is that the set $J$ of homotopy classes of almost complex 
structures on $M$ is equipped with a
structure of an affine space over $H^2(M;\mz)$ or, by Poincar\'e duality, on $H_2(M,\partial M;\mz)$ and one can calculate the
``difference" between two classes $[J_1]$ and $[J_2]$ by
considering the difference between the corresponding
$Spin^c$-structures. In what follows, we study how to explicitly compute this
difference.

Let $M$ and $g$ be as above and let $J_1$ and $J_2$ two orthogonal, positive almost complex
structures on it and $[J_1]$ and $[J_2]$ be their homotopy
classes. 
\begin{prop} \label{comparison}
There exists a well defined class
$\alpha(J_2,J_1):= \alpha([J_2],[J_1])\in H_2(M,\partial M;\mz)$
representing the obstruction to the homotopy between $J_2$ and
$J_1$. The properties of $\alpha$ are:
\begin{enumerate}
\item $\alpha(J_2,J_1)=-\alpha(J_1,J_2)$; \item
$\alpha(\cdot,J_1):J\to H_2(M,\partial M;\mz)$ is bijective for
any $J_1$; \item
$\alpha(J_3,J_2)+\alpha(J_2,J_1)=\alpha(J_3,J_1)$.
\end{enumerate}
Moreover it holds $\alpha(J_2,J_1)=PD(b(J_2)-b(J_1))$ where
$b(J_i),\ i=1,2$ is the $Spin^c$ structure associated to the
homotopy class $[J_i]$ by the canonical bijection $b$ of Theorem \ref{teo:bij}.
\end{prop}

\begin{prf}{1}{
We limit ourselves to sketch how to define $\alpha$. From the
definition, Property $1$ will descend automatically. The other two
properties will descend from Theorem \ref{teo:bij}.

Let us associate to each $J_i$ a self dual $2$-form on $M$ which
$\omega (J_i)$ as follows. At each point $p\in M$ consider a basis of $T_pM$ of the form
$(x_1,J_i(x_1),x_2,J_i(x_2))$ and stipulate that the $2$-form
$\omega(J_i)$ in $p$ is represented by $x_1\wedge
J_i(x_1)+x_2\wedge J_i(x_2)$. The so
obtained $2$-form is a self dual two form whose norm is everywhere
$\sqrt 2$.

First we show that we can choose  $g$-orthogonal $J_1' \in [J_1]$ and
$J_2' \in [J_2]$, so that the set $S=\{p\in M |
J'_1(p)=-J'_2(p)\}=\{p\in M |
\omega_1(p)=-\omega_2(p)\}$ is a properly embedded orientable surface in $(M,\partial M)$.

Since the forms
$\omega_i=\omega(J_i),\ i=1,2$ are self-dual $2$-forms having norm
equal to $\sqrt 2$, they are (up to normalization) sections of the
($2$-dimensional) unit bundle of $\Lambda^+ M$ (the bundle of self
dual two forms) which we will denote by $U\Lambda^+M$. 
The forms $\omega_i$
define embeddings $M_i$ of $M$ into the total
space of a $6$-dimensional bundle. Choosing
$J_i'$ generically we can suppose that the set of points of $M$
where $\omega_2=-\omega_1$ is an orientable surface $(S,\partial
S)$ properly embedded in $(M,\partial M)$. 
Since $M$ is
oriented, also $M_i$ can be oriented by pulling back the
orientation of $M$ through the projection of the fiber bundle
$U\Lambda^+ M$ on $M$. 
Also the total space of the
bundle $U\Lambda^+ M$ can be oriented using the ``horizontal" orientation of $TM$ and the orientation on the fiber fixed by stipulating that, if $e_1,e_2,e_3,e_4$
 is a positive orthonormal basis of $T^*_p M$, the basis of $\Lambda ^+M$ given by $e_1\wedge
e_2 +e_3\wedge e_4, e_1\wedge e_3 -e_2\wedge e_4, e_1\wedge e_4
+e_2\wedge e_3$ is positive.

 Now that all our starting objects are oriented, we
 equip $S$ with an orientation depending only
on the ordered pair $(J_1,J_2)$ as follows. 
Note that $S$ is the  projection of the surface $-M_1\cap M_2$,
 where $-M_1$ is the embedding of $M$ in $\Lambda^+ M$ given by
 $-\omega_1$. Let us orient $-M_1\cap M_2$ locally around a point $q$ as follows. Let $x_1,x_2,x_3,x_4,x_5,x_6$ be a local system of
 coordinates of $U\Lambda^+ M'$ around $q$ such that:
 \begin{enumerate}
 \item the basis $\frac{\partial}{\partial x_1},\frac{\partial}{\partial
 x_2},\frac{\partial}{\partial x_3},\frac{\partial}{\partial x_4},\frac{\partial}{\partial
 x_5},\frac{\partial}{\partial x_6}$ is an oriented basis of the
 tangent space at $q$ of the total space of $U\Lambda^+ M$;
\item the vectors $\frac{\partial}{\partial
x_1},\frac{\partial}{\partial
 x_2},\frac{\partial}{\partial x_3},\frac{\partial}{\partial x_4}$
 form a positive basis of $T_q M_2$;
 \item the vectors $\frac{\partial}{\partial x_3},\frac{\partial}{\partial x_4},\frac{\partial}{\partial
 x_5},\frac{\partial}{\partial x_6}$ form a positive basis of $T_q
 (-M_1)$.
 \end{enumerate}

We stipulate that a positive basis of $T_q (-M_1 \cap M_2)$
is given by $\frac{\partial}{\partial
x_3},\frac{\partial}{\partial x_4}$. This
construction yields a well defined orientation on $-M_1 \cap M_2$
and hence on $S$.

To each component of $S$, we now assign an integer index equal to
$\pm 1$. Let $S_1$ be one component of $S$ and consider a small
disc $D$ transverse to $S_1$ and intersecting it in a single point
$p$; we orient $D$ so to complete the orientation of $S_1$ to the
orientation of $M$. By construction $J_2=-J_1$ on $p$ and not on
$D-\{p\}$. Locally, on $D$, we can identify the section of
$U\Lambda^+M$ given by $J_1$ with the ``north pole'' of the sphere
bundle. Then the image of $D$ through the other section, given by
$-J_2$, is a small disc surrounding the north pole. The index of
$S_1$ is defined to be $1$ is this small disc is oriented as
the spherical fiber and $-1$ otherwise.

Now, we call the class $[S]\in H_2(M,\partial M)$ the {\it
comparison class} between $[J_2]$ and $[J_1]$, and we will denote
it $\alpha([J_2],[J_1])$ or, when it will not cause any confusion,
simply $\alpha(J_2,J_1)$. One can check that $\alpha$ does not depend on the choices made.  Property $1$ descends directly from the
definition. Property $2$ and $3$ are deduced from the
properties of $Spin^c$-structures recalled in Theorem \ref{teo:bij}. They could also be deduced by
independent arguments, as, for instance, the definition of the
analogue of the Pontrjagin-move for almost complex structures and
surfaces in $M$ (see \cite{BP}, Chapter VI).

}
\end{prf}

\subsection{Surjectivity of the reconstruction map}

Next we prove
that all the homotopy classes of almost complex structures can be
obtained via the reconstruction map of Corollary \ref{refinedreconstruction}. 

We first calculate the difference $\alpha(J_{P'},J_P)$ when $P'$ and $P$ are two branched shadows of $M$ connected by a branched move. 
Let us note that we can combinatorially calculate $2\alpha(J_{P'},J_P)$: indeed, by Theorem \ref{teo:bij} and the last statement of Proposition \ref{comparison}, we have $2\alpha(J_{P'},J_P)=c_1(J_{P'})-c_1(J_P)$ and using Corollary \ref{Chern} both $c_1(J_{P'})$ and $c_1(J_P)$ are represented by the sum of the Euler and gleam cochains of $P'$ and $P$ respectively. 
When there is no $2$-torsion in $H^2(M;\mathbb{Z})\cong H_2(M,\partial M;\mathbb{Z})$ we are done since we can divide by two the above equality.
In the following lemma, we calculate the values of $[\Delta Eul]:= 
[Eul(P')]-[Eul(P)]$, $[\Delta gl]:=  [gl(P')]-[gl(P)] \in
H^2(M;\mathbb{Z})$, and
$2\alpha(J_{P'},J_P)\in H_2(M, \partial M;\mathbb{Z})$ when $P'$ is obtained from $P$ by the application of a $0\to 2$-move or a $2\to 3$-move.
\begin{lemma}
If $P'$ is obtained from $P$ through the application of a sliding $0\to 2$-move or $2\to 3$-move (i.e. one of those drawn in the left part of Figure \ref{branched lune} or in the upper part of Figure \ref{branched MP}) then $[\Delta Eul]=[\Delta gl]=2\alpha(J_{P'},J_P)=0$; if the $0\to 2$ or $2\to 3$-move is a bumping one (shown in the right part of Figure \ref{branched lune} or in the lower part of Figure \ref{branched MP}) then $[\Delta Eul]=-2PD(\hat{\Delta})$, $[\Delta gl]=0$ and hence $2\alpha(J_{P'},J_P)=-2PD(\hat{\Delta})$ where $\Delta$ is the upper-right region in our figures of the
bumping moves, the orientation of $\Delta$ is given by the
branching and $PD(\hat{\Delta})$ is the Poincar\'e dual of the element of $H^2(M;\mathbb{Z})$ which is represented by the cochain dual to $\Delta$.
\end{lemma}
\begin{prf}{1}{
The proof is a straightforward computation based on Proposition
\ref{prop:branched presentation} and Corollary \ref{Chern}; as an example, we consider the
calculation for the bumping $2\to 3$-move, shown in the lower part of
Figure \ref{branched MP}. 
Since $P$ and $P'$ describe the same manifold, the presentations of $H^2(M;\mathbb{Z})$ they provide through Proposition \ref{prop:branched presentation} are equivalent, and so, in what follows, we use the presentation provided by $P'$.  

Since the maw is fixed on the boundary of the polyhedron shown in the
figure, its behavior on the boundary of the regions does not change after the move except for $R_1$, $R_2$ and $R_3$. 
For instance, before the move $R_1$ is the preferred
region of the central edge while after the move it is the
preferred region of no edge in the figure. A similar 
phenomenon can be
observed on $\partial R_2$ and $\partial R_3$.
It can be checked that in all these cases the result is the addition
of a full positive or negative twist to the maw on the boundaries of the
regions, so the coefficients of $\hat{R}_1$, $\hat{R}_2$
and $\hat{R}_3$ in the Euler cochain change respectively of $1$, $-1$
and $-1$.
 
Finally, the maw on the boundary of the new region $R$ points always
outwards and hence the index of the singularity on the interior of $R$
is $1$. 
So the difference between the cochains $Eul(P')$ and $Eul(P)$ is
$\hat{R}+\hat{R}_1-\hat{R}_2-\hat{R}_3$; but since the relations in
the presentation of $H^2(M;\mathbb{Z})$ induced by the edges between
$R_1$ and $R_2$ and $R$ and $R_3$ are respectively $\hat{R}_2-\hat{R}_1-\hat{R}_4$ and $\hat{R}_3-\hat{R}-\hat{R}_4$, the total change of the Euler cohomology class is $-2\hat{R}_4=-2\hat{\Delta}$. 
To finish, note that the cohomology class represented by the gleam cochain is not changed by the move because the gleams of the regions involved in the move do not change and the region $R$ has zero gleam.
}\end{prf}

To perform the calculations in the case of $1\to 2$-moves we first set up the notation.
By Lemma \ref{lem:p1branchable}, each $1\to 2$-move is
branchable and has $32$ branched versions, one for each of
the $24$ initial branching of the vertex where the move acts, with 
some cases which split since any orientation of the region
created by the move ($R_7$ in Figure \ref{p1branchable}) gives
a branching. Up to diffeomorphisms of the ball where the
move acts, we reduce to the $16$-cases examined in the proof of Lemma \ref{lem:p1branchable}. Using
the same convention, we summarize in the following table how
do the classes given in $H^2(M;\mathbb{Z})$ by the Euler
cochain, the gleam cochain and the Chern class, change after 
a branched $1\to 2$-move $P\to P'$ in each of the $16$
cases. We use the
presentation of $H^2(M;\mathbb{Z})$ given by $P'$.

\begin{center}\label{table}
\begin{tabular}{|l|l|l|l|l|l|}
\hline Case & $(R_1,R_2,R_3,R_4,R_5,R_6)$ & $R_7$ & [$\Delta$Eul]
& [$\Delta$gl] &
[$\Delta c_1$]\\
\hline
1 &  $(+,+,+,+,+,+)$ & $+$ & 0 & 0 & 0\\
2a & $(+,+,+,+,+,-)$ & $+$ & $-\hat{R}_6$ & $-\hat{R}_6$ & $-2\hat{R}_6$ \\
2b & $(+,+,+,+,+,-)$ & $-$ & $-\hat{R}_5$ & $-\hat{R}_5$ & $-2\hat{R}_5$ \\
3a & $(+,+,+,+,-,+)$ & $+$ & $-\hat{R}_5$ & $-\hat{R}_5$ & $-2\hat{R}_5$ \\
3b & $(+,+,+,+,-,+)$ & $-$ & $-\hat{R}_6$ & $-\hat{R}_6$ & $-2\hat{R}_6$ \\
4 &  $(+,+,+,+,-,-)$ & $-$ & 0 & 0 & 0\\
5 &  $(+,+,+,-,+,+)$ & $+$ & $\hat{R}_4$ & $-\hat{R}_4$ & 0\\
6 &  $(+,+,-,-,+,+)$ & $+$ & $\hat{R}_5$ & $-\hat{R}_5$ & 0\\
7a & $(+,+,-,+,-,+)$ & $+$ & $-\hat{R}_4$ & $-\hat{R}_4$ & $-2\hat{R}_4$\\
7b & $(+,+,-,+,-,+)$ & $-$ & $-\hat{R}_2$ & $-\hat{R}_2$ & $-2\hat{R}_2$\\
8 &  $(+,-,+,+,-,-)$ & $-$ & $\hat{R}_2$ & $-\hat{R}_2$ & 0\\
9 &  $(+,+,-,-,-,+)$ & $+$ & 0 & 0 & 0\\
10 & $(+,-,-,+,-,-)$ & $-$ & $\hat{R}_6$ & $-\hat{R}_6$ & 0\\
11 & $(+,-,-,+,-,+)$ & $-$ & 0 & 0 & 0\\
12a & $(+,-,-,-,-,+)$ & $+$ & $-\hat{R}_2$ & $-\hat{R}_2$ & $-2\hat{R}_2$\\
12b & $(+,-,-,-,-,+)$ & $-$ & $-\hat{R}_4$ & $-\hat{R}_4$ & $-2\hat{R}_4$\\
\hline
\end{tabular}
\end{center}

The above computations give the value of
$2\alpha(J_{P'},J_P)=[\Delta c_1]:=  c_1(J_{P'})-c_1(J_P)$; note that this value has always the form $2\hat{R}$ for a suitable region $R$ and so $[\hat{R}]$ is a natural candidate for the class $\alpha(J_{P'},J_P)$. 
\begin{teo}\label{valoredialfa}
Let $P$ and $P'$ be two branched shadows of the same manifold
connected by a $0\to 2$, $2\to 3$ or $1\to 2$ branched move, and let $c_1(J_{P'})-c_1(J_P)=2\hat{R}$ for a suitable region $R$ be the difference of the Chern classes of the associated almost complex structures, calculated as explained above. Then $\alpha(J_{P'},J_P)=R$.
\end{teo}
\begin{prf}{1}{
The application of a branched move changes a shadow only in a ball $B^4$ of $M$ and hence the class $\alpha(J_{P'},J_P)$ is represented by a properly embedded oriented surface $S$ contained in $B^4$ whose connected components are equipped with multiplicities (see Proposition \ref{comparison}). Up to a small isotopy, its intersection with $P'$ (or $P$) can be supposed to be made of a finite number $n_i$ of points in the interior of each region $R_i$; $PD(\alpha(J_{P'},J_P))=PD([S])=\Sigma_i n_i\hat{R}_i$,  hence the Poincar\'e dual of $\alpha(J_{P'},J_P)$ can be represented in $H^2(M;\mathbb{Z})$ by a cochain whose coefficients are non zero only on the regions intersecting $B^4$.
 
We claim that there is no class $\beta\in H^2(M;\mathbb{Z})$ other than $\hat{R}$ such that $2\beta=2\hat{R}$, whose coefficients are zero on the regions not intersecting $B^4$. Suppose the contrary, then $2(\beta-\hat{R})=0$ and  we would have a $2$-cochain $\gamma$ representing an element of $2$-torsion in $H^2(M;\mathbb{Z})$  whose coefficients are zero on the regions not intersecting $B$. Moreover the relations used to show that $2\gamma=0$ should all come from edges intersecting $B^4$ since otherwise we would use informations on $M$ not inherent to the move but also on other parts of $M$ not involved in the move itself. Let $Z$ be the $\mathbb{Z}$-module generated by the cochains $\hat{R}_i$ subject to the relations induced by the edges $e_j$ where $i$ and $j$ range respectively over the regions and the edges of $P$ intersecting $B^4$. All the above cochains ($\gamma$ included) are naturally represented in $Z$ and in particular it should hold $2\gamma=0$ in $Z$, but $Z$ has no $2$-torsion, and we are done. 

We show here that $Z$ has no $2$-torsion when $P'$
and $P$ are connected by a $1\to 2$-move and leave the cases of the $2\to 3$ and $0\to 2$ moves to the reader. Using the notation of the left part of Figure \ref{p1branchable}, the module $Z$ is generated by $\hat{R}_i,\ i=1,..,6$ and has $4$-relations namely: $\hat{R}_2=\hat{R}_1+\hat{R}_5$, $\hat{R}_3=\hat{R}_5+\hat{R}_4$, $\hat{R}_2=\hat{R}_6+\hat{R}_3$, $\hat{R}_1=\hat{R}_6+\hat{R}_4$. Then $Z$ is isomorphic to $\mathbb{Z}\oplus\mathbb{Z}\oplus\mathbb{Z}$.

}\end{prf}

\begin{teo}\label{teo:surjectivityofrefinedreconstruction}
The refined reconstruction map from branched shadows of $M$ to pairs
$(M,[J])$ with $[J]$ homotopy class of positive almost complex structures on $M$, is surjective.  
\end{teo}
\begin{prf}{1}{
Fix an auxiliary riemannian metric $g$ on $M$. We limit ourselves to give an idea of the proof since it is an adaptation of that of Theorem 4.6.4 of \cite{BP}. Let $[J]$ be as in the statement, $P$ be a branched shadow of $M$ (which exists by Proposition
\ref{prop:existsbranchedshadow}) and let $[J_P]$ the homotopy class associated to $P$ by Corollary \ref{refinedreconstruction}. We want to show that there exists a branched shadow $P'$ such that $[J_{P'}]=[J]$; if $[J]=[J_P]$ we are done. 
Otherwise, let us write the Poincar\'e dual of $\alpha(J,J_P)$ as $\sum_i
k_i \hat{R}_i$ where $\hat{R}_i$ is the cochain dual to the region
$R_i$ of $P$ (using the presentation of $H^2(M;\mz)$
given by $P$). The idea of the proof is to apply a suitable
sequence of $0\to 2$ and $2\to 3$-moves to $P$ to modify it and get a new branched shadow of $M$ carrying
the homotopy class $[J]$. For that, it is
sufficient to exhibit a sequence that ``decreases" the difference
$\alpha$ between the $J$ and $J_P$ by a cochain whose Poincar\'e dual is cohomologous to one of the form $\pm
\hat{R}_i$ for any region $R_i$ of $P$. 
For instance, suppose that
we want to get a shadow carrying a homotopy class which differs from $J_P$
by $-\hat{R}_i$, and that the boundary of the region
$R_i$ contains an edge $e_j$ whose preferred region is not $R_i$; then,
by Theorem \ref{valoredialfa} a self $0\to 2$-move is a bumping
move producing a branched shadow $P'$ such that
$\alpha(J_{P'},J_P)=-\hat{R}_i$. 
If for all the edges in $\partial
R_i$ the preferred region is $R_i$, let $e_j$ be an edge of
$\partial R_i$ and $R_k$ and $R_l$ be other regions containing
it. Since $\hat{R}_i=\hat{R}_k+\hat{R_l}$ it is sufficient to
apply the above moves both to $R_k$ and $R_l$ to get the wanted
shadow. The sequence of moves producing a difference of the form
$+\hat{R}_i$ is more complicated and we refer to \cite{BP} for a complete account.
}
\end{prf}
\section{Branched shadows and complex structures}
\begin{defi}
An almost complex structure $J$ on a smooth $4$-manifold $M$ is said
to be {\it integrable} or {\it complex} if for each point $p$ of $M$
there is a local chart of $M$ with values in $\mathbb{C}^2$ transforming $J$ into the complex structure of $\mathbb{C}^2$.
\end{defi}
\subsection{Branched shadows in complex manifolds}
In this subsection, supposing that $M$ is equipped with an integrable structure $J$,
we adapt to the case of branched shadows a series of classical results of Bishop \cite{Bi},
Chern and Spanier \cite{CS}, and  Harlamov and Eliashberg
\cite{HE} regarding invariants of embeddings of real surfaces in
complex manifolds.

Let $P$ be a branched shadow embedded in $M$. Up to perturbing the embedding of $P$ through a small isotopy we can suppose that there is only a finite
number of points $p_1,\ldots, p_n$ and $q_1,\ldots,q_m$, contained in the regions of $P$ where $T_{p_i}P$ (resp. $T_{q_j}P$) is a complex plane such that the orientations induced by the branching of $P$ and by the complex structure coincide (resp. do not coincide).

\begin{defi}
The points $p_1,\ldots,p_n$ are called
{\it positive complex} points of $P$ or simply {\it positive}
points. Analogously, the points $q_1,\ldots,q_m$ are called {\it
negative complex} points of $P$ or {\it negative} points. All the
other points of $P$ are called {\it totally real}.
\end{defi}

To each complex point $p$ of a region $R_i$ of $P$ we can assign an integer number
called its {\it index}, denoted $i(p)$, as follows. Fix a small disc $D$ in $R_i$ containing $p$ and no
other complex point and let $N$ be the radial vector field
around $p$. The field $J(N)$ on
$\partial D$ is a vector field transverse to $P$ since no point on
$D-p$ is complex. Let $\pi (J(N))$ be the projection of this field
onto the normal bundle of $D$ in $M$. Since $D$ is contractile,
this bundle is trivial and we can count the number $\nu(p)$
of twists performed by $\pi(J(N))$ while following $\partial D$
($D$ and $M$ are oriented). The index of $p$ is: $i(p)=\nu(p)+1$. 
Moreover, we define $\nu(R_i)$ as the sum
over all the complex points $p$ of $R_i$ of $\nu(p)$. 

Up to a small perturbation by an isotopy of the embedding of $P$ in $M$ we can assume that all the indices of the complex points of $P$ are equal to $\pm 1$.

\begin{defi}
A complex point $p$ of $P$ whose index is equal to $1$ is {\it elliptic}, if its index is $-1$ it is {\it hyperbolic}.
\end{defi}

We define the index $c_1(p)$ associated to each complex point $p$ of a region $R_i$ of $P$ as follows. 
Let $D$ and $N$ be as above and complete $N$ on $\partial D$ to a basis of $TD$ by using
the field $T=T\partial D$ tangent to the boundary of $D$. The pair
of fields $(N,T)$ gives a basis of $TD$ in each point $q$
of $\partial D$, and, since no such point is complex, they can
be completed to a positive complex basis of $T_qM$ given by 
$(N,J(N),T,J(T))$. Let now $\frac{\partial}{\partial z}$ and
$\frac{\partial}{\partial w}$ be two vector fields defined on a
neighborhood of $D$ in $M$ such that $(\frac{\partial}{\partial
z},\frac{\partial}{\partial w})$ is pointwise a complex basis of
$TM$. Then, on each point $q$ of $\partial D$ we can compare the
two complex bases given by $(N+J(N),T+J(T))$ and
$(\frac{\partial}{\partial z},\frac{\partial}{\partial w})$ by
considering the determinant $det_q$ of the change of basis from
the latter to the former basis. The value of the index of $det_q$ around $0$ in
$\mathbb{C}$ while $q$ runs across $\partial D$ according to the
orientation of $D$, is defined to be $c_1(p)$. We define $c_1(R_i)$ as the sum of $c_1(p)$ over all the complex points $p$ of $R_i$.

The following result due to Bishop (\cite{Bi}), Chern and Spanier
(\cite{CS}) and Lai \cite{Lai}, shows how the above defined indices
are related to the topology of $R_i$ and $M$:
\begin{teo} \label{teo:bishop}
Let $R_i$ be a surface with boundary contained in complex manifold $M$ such that $\partial R_i$ does not contain complex points and let $I^+=\sum_{i=1,\ldots,n}i(p_i)$ and
$I^-=\sum_{j=1,\ldots,m}i(q_j)$. Then the following equalities
hold:
$$I^+=\frac{1}{2}(\chi (R_i)+\nu(R_i)+c_1(R_i))$$
$$I^-=\frac{1}{2}(\chi(R_i)+\nu(R_i)-c_1(R_i))$$
\end{teo}

We will need the following, due to Harlamov and
Eliashberg \cite{HE}:
\begin{teo}[Annihilation theorem]\label{teo:harlamov-eliashberg}
Let $S$ be an oriented real surface embedded in a complex manifold $M$ and let $p_1$ and $p_2$ be two complex
points of $S$ of the same sign (i.e. both positive or negative) and belonging to the same connected component of $S$.
Let $\alpha$ be an arc in $S$ connecting $p_1$ and $p_2$ and
containing no other complex point of $S$ and suppose that
$i(p_1)=1=-i(p_2)$. 
There exists on $S$
a small isotopy $\phi_t,\ t\in [0,1]$ which is the identity out of a small neighborhood $U(\alpha)$ of
$\alpha$ and such that $\phi_1(U(\alpha))$ contains no complex points.
\end{teo}

The following is the analogous in the world of shadows of Theorem \ref{teo:harlamov-eliashberg}:
\begin{lemma}\label{lem:branchedHE}
Let $R_i$, $R_j$ and $R_k$ be three regions of $P$ adjacent along
a common edge $e\in Sing(P)$ so that $R_i$ is the preferred region of $e$. 
There exists an isotopy
$\phi_t:P\to M,\ t\in [0,1]$ whose support is contained in a small
ball $B$ around the center of $e$ such that $\phi_1(B\cap P)$ contains three
more complex points $p_i$, $p_j$ and $p_k$ respectively in $R_i$,
$R_j$ and $R_k$ whose indices are respectively $\pm 1$, $\mp 1
$ and $\mp 1$.
\end{lemma}
\begin{prf}{1}{
By modifying with a $C^0$-small isotopy the position of $R_i$ in a
neighborhood of an interior point we can assume that, near the
edge $e$ there is a pair of complex points $p_i$ and
$p'_i$ in $R_i$ of opposite indices (first create some complex points through a small isotopy and then apply Harlamov-Eliashberg result to delete all of them but two). Roughly speaking, we now ``slide $R_j$ over'' $p_i'$ and show that both on
$R_k$ and $R_j$ two complex points are created by this sliding.

Consider an arc $\gamma$ embedded in $R_i$
whose endpoints are $p'_i$ and a point $q$ on $e$. Let $D$ be a small disc neighborhood
in $R_j$ of $q$ and let $D'$ a slightly bigger neighborhood of
$D\cup\beta$ in $R_i \cup R_j$ (since $P$ is branched,
$D'$ is a smooth disc).  Consider the isotopy that fixes every
point of $P-D$ and moves $D$ by letting it slither over $\beta$
and pass over $p'_i$. This isotopy acts in a $4$-ball and changes only the regions $R_i$ and $R_k$ along their borderlines. The point $p'_i$ passes from 
one side to the other one, becoming a
complex point $p_k$ in $R_k$;
we are left to prove that a $J$-complex point has been
created on $R_j$ having the same index as $p_k$.
This is proved by using Theorem \ref{teo:bishop}. 
Consider the image $D''$ of the bigger disc $D'$ after the
isotopy: it is a disc whose boundary is made only of totally real
points and $i(D'')=i(D')$ since $\partial D''=\partial D'$ by
construction. Then we finish by observing that $i(D')=i(p'_i)$
since $D'$ by construction contains only $p'_i$.
}\end{prf}

The above lemma suggests the following:
\begin{defi}
The {\it positive index} and {\it negative index} cochains of $P$, denoted respectively $I^+(P)$ and $I^-(P)$ are the $2$-cochains given by $\Sigma_i I^{\pm}(R_i)\hat{R}_i$, where $i$ ranges over all the regions of $P$.
\end{defi}
\begin{teo}
The cohomology classes $[I^{\pm}(P)]\in H^2(M;\mathbb{Z})$ are invariants of the embedding of $P$ in $M$ up to isotopy.
\end{teo}
\begin{prf}{1}{
Let $\phi_t, t\in [0,1]$ be an isotopy of $P$ in $M$ so that $\phi_0=id$ and $\phi_1(P)=P'$. Up to slightly perturbing $\phi$, we can suppose that the following holds:
\begin{enumerate}
\item the number of creation/annihilations of complex points with opposite indices during the isotopy is finite;
\item no creation/annihilation of complex points at time $t$  happens on $\phi_t(Sing(P)), t\in [0,1]$;
\item the complex points of $\phi_t(P)$ cross $\phi_t(Sing(P))$ only a finite number of times and transversally in the interior of the edges of $\phi_t(Sing(P))$.
\end{enumerate}
Then we have to check invariance of the classes $I^{\pm}(P)$ under two kinds of catastrophes: when a pair of complex points is created or annihilated in the interior of a region of $P$ and when a complex point crosses $\phi_t (Sing(P))$. Invariance in the first case comes from Theorem \ref{teo:harlamov-eliashberg} and the definition of $I^{\pm}(R_i)$; the second case is a consequence of Lemma \ref{lem:branchedHE}.
}\end{prf}

We now compare the almost
complex structure carried by a branched shadow $P$ with the ambient complex structure $J$.
\begin{prop}\label{prop:spostapunti}
The following holds: $\alpha(J,J_P)=PD([I^-(P)])$ where $PD$ is the
Poincar\'e dual of the cohomology class represented by the negative
index cochain of $P$. Moreover, if $\alpha(J,J_P)=0$ there exists an isotopy $\phi_t,\ t\in[0,1]$ of $P$ in $M$ such
that $\phi_0$ is the identity and $\phi_1(P)$ contains no
negative $J$-complex points. 
\end{prop}
\begin{prf}{1}{

Up to homotopy, we can suppose that $J_P$ is generic with respect to $J$ and 
represent $\alpha(J,J_P)$ as an oriented properly embedded surface
$S$ in $(M,\partial M)$ whose components are equipped with indices
equal to $\pm 1$. Consider the $2$-cochain $\beta= \sum n_i \hat{R}_i$ where $i$ ranges
over all the regions $R_i$ of $P$ and $n_i$ is equal to the sum of
the indices of all the intersection points between $S$ and $R_i$.
By the construction of $S$ (see Proposition \ref{comparison}), 
$P\cap S$ is the set of negative $J$-complex points of $P$ equipped with indices equal to those of the
corresponding components of $S$, hence $\beta=I^-(P)$. 

On the other side, the evaluation of $\beta$ on a cycle $C\in H_2(M;\mathbb{Z})$
is equal to the intersection of $C$ and $S$ in $M$ and hence $\beta=PD(\alpha(J,J_P))$: indeed $C$ can
be represented by a surface lying in a small neighborhood
of $P$ in $M$ and intersecting $S$ only near $P\cap S$. 

If $\alpha(J,J_P)=0$, then $[I^-(P)]=0$ in $H^2(M;\mathbb{Z})$
and its expression as a cochain can be reduced to $0$ by
means of the relations given by the edges of $P$ (see
Proposition \ref{prop:branched presentation}). To conclude it is
sufficient to prove that, given an arbitrary edge $e$ of $Sing(P)$
inducing a relation of the form $\hat{R}_i=\hat{R}_j+\hat{R}_k$ on
the three regions containing it, it is possible to find an isotopy
$\phi_e$ of $P$ in $M$ which has the effect of eliminating a point
$p'_i$ of index $\pm 1$ on $R_i$ and creating a pair of points of
the same index $p_k$ and $p_j$, respectively on $R_k$ and on
$R_j$. This is the statement of Lemma
\ref{lem:branchedHE}.
 }
\end{prf}
\subsection{Integrable representatives of almost complex structures} 
This subsection is devoted to prove Theorem \ref{mainteo}:

{\it Proof of {\rm \ref{mainteo}}.}
Let $[J]$ be a homotopy class of almost complex structures on $M$. By Theorem \ref{teo:surjectivityofrefinedreconstruction} there exists a branched shadow $P$ of $M$ such that the homotopy class of the almost complex structure $J_P$ associated to $P$ is $[J]$. 
To construct an integrable representative of $[J]$ we now immerse
$M$ in a suitable complex manifold and pull back its complex structure
to $M$ through the immersion (the pull-back through a local
diffeomorphism of a complex structure is a complex structure). We then prove that
the so obtained complex structure belongs to $[J]$. The branched
shadow $P$ is crucial both in the construction of the immersion
of $M$ and in the control of the homotopy class of the complex structure.

First, since $Sing(P)$ is a graph, its regular
neighborhood $N$ in $P$ can be embedded in $\mathbb{C}^2$ equipped
with coordinates $(z,w)$. Moreover, we can suppose that the projection
over the plane $w=0$ has surjective  and positive differential (i.e. the image of the orientation of $TP$ coincides with that of the plane $w=0$) in every point of $N$. This implies that the image of $N$ in $\mathbb{C}^2$ contains no negative complex points. 
Let us now extend arbitrarily the so constructed immersion of $N$ in $\mathbb{C}^2$ to an immersion of the whole $P$: this can be done since each boundary component of $N$ coincides with a boundary component of a region of $P$ and each loop in $\mathbb{C}^2$ bounds an immersed disc. 

That way we construct an immersion of $P$ in $\mathbb{C}^2$ but we have to solve two problems:
\begin{enumerate}
\item modify the immersion so to delete all the negative complex points of the image;
\item extend the new immersion to the whole $M$. 
\end{enumerate}

We solve both problems with the same technique, but to use it we
first fix the notation. Up to isotopy of the immersion we can
suppose that for each region $R_i$ of $P$, the image of $R_i$ contains
two small discs $D_i^+$ and $D_i^-$ parallel to the $w=0$ plane and oriented by $R_i$ respectively in the positive and in the negative way with respect to the complex orientation of the plane. Let $p_i^+$ and $p_i^-$ points in $D_i^+$ and $D_i^-$ respectively and let $C'$ be the complex $4$-manifold obtained by deleting small neighborhoods of $p_i^{\pm}$ of the form $p_i^{\pm}+\{(z,w)| ||z||\leq \epsilon_i, ||w||\leq \epsilon_i\}$ from $\mathbb{C}^2$. 
If we glue back to $C'$ along $\partial D_i^-$ a $2$-handle equipped with a complex structure of the form $\{(z,w)|\  ||z||\leq \epsilon_i ,||w||\leq\frac{\epsilon_i}{2}\}$ though the map $(z,w)\rightarrow p_i^-+(z,z^{n_i}w)$ and acting analogously for $D_i^+$ (using an exponent $m_i$), we obtain an immersion of $P$ in a new complex manifold by sending $D_i^{\pm}$ to the cores of the $2$-handles. If we slightly perturb this new immersion we see that the points $p_i^{\pm}$ are isolated complex points (respectively of positive and negative type)  whose indices are $m_i$ and $n_i$. 
The above recipe allows us to modify the immersion in small neighborhood of points of each region (changing also the codomain of the immersion) and change arbitrarily the total indices $I^\pm$ of the image of each region. Then, we can fix $I^-(R_i)$ to be zero on each region, and, after applying a suitable number of times Theorem \ref{teo:harlamov-eliashberg}, we can suppose that each $R_i$ contains no negative points.

Let us now solve the second problem; the pull back of a neighborhood of the image of $P$ through the so obtained immersion gives a thickening of $P$ but a priori it is not true that such a thickening coincides with $M$. It is so if the gleams induced on the regions of $P$ by this thickening (recall the definition of the gleam) are the same as those induced by $M$ on $P$. Changing suitably the coefficients $m_i$ in the above construction we can change arbitrarily the gleams induced on each region $R_i$ without adding negative complex points; indeed the gleam is a relative version of the self-intersection number of each region and adding a twist to the attaching map of the $2$-handle around $p_i^+$ changes by one this self intersection number.    

Then, up to suitably reattaching to $C'$ all the $2$-handles corresponding to the points $p_i^{\pm}$, the immersion of $P$ can be extended to an immersion of $M$ in $C''$ and the pull back of the complex structure of $C''$ to $M$ is a complex structure such that $P$ contains no negative point.
By Proposition \ref{prop:spostapunti} we conclude that the so obtained complex structure on $M$ belongs to the homotopy class carried by $P$ and we are done. 
 \setcounter{teo}{2}
 \setcounter{section}{1}
\hfill\qed

\end{document}